\documentclass[12pt]{article}
\usepackage{mycommands}
\usepackage{colonequals}
\usepackage{tikz}
\usetikzlibrary{decorations.pathreplacing}
\title{Degenerations of Curves in Projective Space and the Maximal Rank Conjecture}

\begin{document}
\maketitle

\begin{abstract}
In this note, we give an overview of a new technique for studying
Brill--Noether
curves in projective space via degeneration. In particular,
we give a roadmap to 
the proof of the Maximal Rank Conjecture.     
\end{abstract}

\section{Introduction}

The technique of degeneration to a reducible
curve has enabled the proof of many results in the theory of algebraic curves.
These results include the \emph{Brill--Noether theorem},
proven by Griffiths and Harris \cite{bn}, Gieseker \cite{gp}, Kleiman and Laksov \cite{kl}, and others,
which describes the space of maps from a general curve to projective space:
If $C$ is a general curve of genus $g$, 
it states that there exists a nondegenerate degree $d$ map $C \to \pp^r$ if and only if
the \emph{Brill--Noether number $\rho(d, g, r)$} is nonnegative, where
\[\rho(d, g, r) \colonequals (r + 1) d - rg - r(r + 1).\]
Moreover, in this case, there exists a unique component $\bar{M}_g^\circ(\pp^r, d)$
of Kontsevich's space of stable maps $\bar{M}_g(\pp^r, d)$
that both
dominates the moduli space of curves
$\bar{M}_g$ and whose general member is nondegenerate.
We call curves in this component of the
space of stable maps
\emph{Brill--Noether curves (BN-curves)}.

These results have been extended in various ways. For example,
Sernesi argues by degeneration to produce components of $\bar{M}_g(\pp^r, d)$ whose image in $\bar{M}_g$
is of the expected dimension when the Brill--Noether number is negative \cite{sernesi}.

However, in this note we focus on the question: \emph{How can we study the geometry
in projective space
of the general BN-curve via degeneration?}
Our goal here is to give an overview of a series of papers by the author and others
\cite{aly, quadrics, ibe, tan, rbn, rbn2, p4, vogt} which give rise to
a general technique for studying BN-curves via degeneration.
This technique is then applied in \cite{mrc} (in conjunction with results
on hyperplane sections of BN-curves obtained in \cite{hyp})
to give a proof of the \emph{Maximal Rank Conjecture},
a conjecture made originally by Severi in 1915 \cite{severi} which determines the Hilbert function
of a general BN-curve:

\begin{conj}[Maximal Rank Conjecture] \label{c:mrc}
If $C \subset \pp^r$ is a general BN-curve ($r \geq 3$), the restriction maps
\[H^0(\oo_{\pp^r}(k)) \to H^0(\oo_C(k))\]
are of maximal rank (i.e.\ either injective or surjective).

Equivalently, the dimension of the space of polynomials
of degree $k$ which vanish on $C$ is given by
\[\begin{cases}
\binom{r + k}{k} - (kd + 1 - g) & \text{if $kd + 1 - g \leq \binom{r + k}{k}$ and $k \geq 2$;} \\
0 & \text{otherwise.}
\end{cases}\]
\end{conj}
The Maximal Rank Conjecture can also be reformulated cohomologically: From the long exact
sequence in cohomology arising from the short exact sequence of sheaves
\[0 \to \ii_{C \subset \pp^r} (k) \to \oo_{\pp^r}(k) \to \oo_C(k) \to 0,\]
we see that $C$ satisfies ``maximal rank for polynomials of degree $k$''
if and only if
\[H^0(\ii_{C \subset \pp^r}(k)) = 0 \tor H^1(\ii_{C \subset \pp^r}(k)) = 0.\]

Many special cases of the maximal rank conjecture have been previously studied,
using an approach originally due to Hirschowitz:
Degeneration to a reducible curve $C' \cup C''$ with $C''$ contained in a hypersurface $S$ of degree $n$
(typically a quadric if $r = 3$ and a hyperplane if $r \geq 4$),
and $C'$ transverse to $S$:

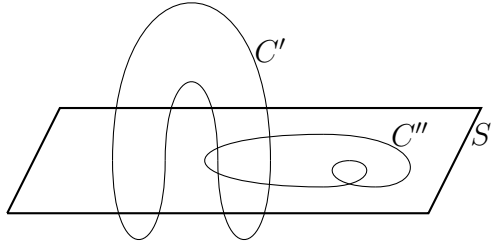
\begin{center}
\begin{tikzpicture}[scale=0.7]
\draw[thick] (0, 1) -- (1, 3) -- (9, 3) -- (8, 1) -- (0, 1);
\draw (2, 2) .. controls (2, 6) and (5, 6) .. (5, 2);
\draw (2, 2) .. controls (2, 0) and (3, 0) .. (3, 2);
\draw (4, 2) .. controls (4, 0) and (5, 0) .. (5, 2);
\draw (3, 2) .. controls (3, 4) and (4, 4) .. (4, 2);
\draw (6, 2.5) .. controls (3, 2.5) and (3, 1.5) .. (6, 1.5);
\draw (6, 2.5) .. controls (8, 2.5) and (8, 1.5) .. (7, 1.5);
\draw (6, 1.5) .. controls (7, 1.5) and (7, 2) .. (6.5, 2);
\draw (7, 1.5) .. controls (6, 1.5) and (6, 2) .. (6.5, 2);
\draw (9, 2.5) node{$S$};
\draw (5, 4.1) node{$C'$};
\draw (7.65, 2.5) node{$C''$};
\end{tikzpicture}
\end{center}

In this case, from the long exact sequence in cohomology arising
from the short exact sequence of sheaves
\[0 \to \ii_{C' \subset \pp^r}(k - n) \to \ii_{C' \cup C'' \subset \pp^r}(k) \to \ii_{C'' \cup (C' \cap S) \subset S}(k) \to 0,\]
we conclude that to show $H^i(\ii_{C' \cup C'' \subset \pp^r}(k)) = 0$ as desired,
it suffices to show
\[H^i(\ii_{C' \subset \pp^r}(k - n)) = H^i(\ii_{C'' \cup (C' \cap S) \subset S}(k)) = 0.\]
One can thus hope to argue by induction on $r$ (if $S \simeq \pp^{r - 1}$ is a hyperplane) and $k$.
However, three fundamental difficulties have limited this approach to special cases:
\begin{enumerate}
\item We need a uniform way to construct the degenerations $C' \cup C''$.
Previous methods were ingenious, but relatively ad-hoc, and hence
not generalizable.

\item It is not possible to always find such reducible curves at which
the fiber dimension of the map $\bar{M}_g(\pp^r, d) \to \bar{M}_g$ at $[C' \cup C'']$
is $\rho(d, g, r) + \dim \aut \pp^r$. We therefore need some other way to see
that such curves are BN-curves.

\item This approach relates maximal rank for $C' \cup C''$ to maximal rank for $C'$
and maximal rank for $C'' \cup (C' \cap S)$. But 
$C'' \cup (C' \cap S)$ is not a curve, so we need a stronger inductive hypothesis.

Even worse, $C'$ and $C''$ must satisfy various incidence
conditions, so $C''$ and $C' \cap S$ are not \emph{independently}
general and there is no nice description of $C'' \cup (C' \cap S)$ that doesn't reference
the \emph{entire} reducible curve $C' \cup C''$.
\end{enumerate}

We begin by discussing the first two difficulties (in Sections~\ref{sec:uniform} and~\ref{sec:fiber}),
which arise whenever we wish to study the geometry of general 
BN-curves via degeneration.
Namely, we show the existence of such degenerations of BN-curves can be
reduced to the existence of integer solutions to certain systems of inequalities.

Then we discuss the third difficulty (in Section~\ref{sec:hyp}),
which is specific to the proof of the Maximal Rank Conjecture.

Finally (in Section~\ref{sec:inequalities}),
we describe a method for proving the existence of integer solutions
to the type of systems of inequalities that appear when applying this method
to the maximal rank conjecture.

\section{The Uniform Construction of Reducible Curves \label{sec:uniform}}

We construct our desired reducible curves via the following method.
First, we fix a finite set of points $\Gamma$ which is general in $\pp^r$,
or general in a hyperplane (or other hypersurface of small degree)
$H \subset \pp^r$. Then, we find BN-curves $C' \subset \pp^r$, and $C'' \subset \pp^r$;
or $C' \subset \pp^r$ transverse to $H$, and $C'' \subset H$ ---
both passing through $\Gamma$:
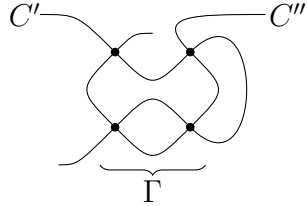
\begin{center}
\begin{tikzpicture}
\filldraw (1, 1) circle [radius=0.05];
\filldraw (2, 1) circle [radius=0.05];
\filldraw (1, 2) circle [radius=0.05];
\filldraw (2, 2) circle [radius=0.05];
\draw (0, 2.5) .. controls (0.5, 2.5) .. (1, 2);
\draw (0.25, 0.5) .. controls (0.5, 0.5) .. (1, 1);
\draw (1, 1) .. controls (1.5, 1.5) .. (2, 1);
\draw (1, 2) .. controls (1.5, 1.5) .. (2, 2);
\draw (2, 1) .. controls (3, 0) and (3, 3) .. (2, 2);
\draw (3, 2.5) .. controls (2, 2.5) and (1.5, 2.5) .. (2, 2);
\draw (2, 2) .. controls (2.5, 1.5) .. (2, 1);
\draw (2, 1) .. controls (1.5, 0.5) .. (1, 1);
\draw (1, 1) .. controls (0.5, 1.5) .. (1, 2);
\draw (1, 2) .. controls (1.25, 2.25) .. (1.5, 2.25);
\draw (-0.2, 2.5) node{$C'$};
\draw (3.3, 2.5) node{$C''$};
\draw [decorate, decoration={brace, mirror, amplitude=0.75ex}] (0.8, 0.5) -- (2.2, 0.5);
\draw (1.5, 0.15) node{$\Gamma$};
\end{tikzpicture}
\end{center}
Taking their union
then gives a reducible curve $C = C' \cup C''$ as desired.

To carry out this construction, we seem to need an answer
to the questions: \emph{When does there exist a BN-curve of given degree $d$
and genus $g$ passing through a set of $n$ general points in $\pp^r$?
If some of these points are constrained to lie in a hypersurface of small degree
(usually a hyperplane), can we find such a BN-curve
transverse to this hypersurface?}
We will not be able to answer these questions, but we will get close enough for our needs.

For this first question, we are asking when the natural map $\pi \colon \bar{M}_{g,n}^\circ(\pp^r, d) \to (\pp^r)^n$
is dominant. The natural conjecture is that $\pi$ is dominant if and only if the dimensions allow it:
\[(r + 1)d - (r - 3)(g - 1) + n = \dim \bar{M}_{g,n}^\circ(\pp^r, d) \geq \dim (\pp^r)^n = rn,\]
or upon rearrangement, if and only if
\begin{equation} \label{nbound}
n \leq \left\lfloor \frac{(r + 1)d - (r - 3)(g - 1)}{r - 1}\right \rfloor.
\end{equation}

This condition can also be expressed as
$\chi(N_C(-p_1 - p_2 - \cdots - p_n)) \geq 0$, where $N_C$ denotes
the normal bundle of $C$, and $p_1, p_2, \ldots, p_n \in C$ are the $n$
marked points.
Since the obstruction to smoothness of $\pi$ lies in $H^1(N_C(-p_1 - p_2 - \cdots - p_n))$, this follows in turn from
the following property 
for the normal bundle $N_C$:

\begin{defi}
We say that a vector bundle $\mathcal{E}$ on an irreducible curve $C$
\emph{satisfies interpolation} if for a general effective
Cartier divisor
$D \subset C$ of every nonnegative degree, either
\[H^0(\mathcal{E}(-D)) = 0 \tor H^1(\mathcal{E}(-D)) = 0.\]
\end{defi}

The property of interpolation for normal bundles is studied
in the following sequence of papers:

\begin{enumerate}[\bf A.]
\item \label{d:aly} In joint work with Atanasov and Yang \cite{aly}, we show that the normal bundle of a general
nonspecial BN-curve (i.e.\ one with $d \geq g + r$) satisfies interpolation
except in exactly three cases: $(d, g, r) \in \{(5, 2, 3), (6, 2, 4), (7, 2, 5)\}$.
Even though the normal bundle of a general BN-curve of degree $6$ and genus $2$ in $\mathbb{P}^4$
does not satisfy interpolation,
it turns out that such curves can still pass through the expected number of general points.
We conclude that a general nonspecial
BN-curve passes through $n$ general points if and only if
\eqref{nbound} holds, with exactly two exceptions: $(d, g, r) \in \{(5, 2, 3), (7, 2, 5)\}$. 

The argument is by inductive degeneration of $C$ to a reducible nodal curve
$X \cup Y$. Older results of Hartshorne and Hirschowitz \cite{hh}
give geometric descriptions of $N_{X \cup Y}|_X$ and $N_{X \cup Y}|_Y$;
however, to describe $N_{X \cup Y}$, one needs a compatible description of
the gluing data $N_{X \cup Y}|_X|_{X \cap Y} \simeq N_{X \cup Y}|_Y|_{X \cap Y}$,
which is quite difficult in general.

The key new insight is to study line subbundles of the normal bundle
obtained by saturating the images of vertical tangent spaces of projection maps.
These enable
us to give an essentially complete geometric description of the gluing data
when $Y$ is a line.

\item \label{d:quadrics} In \cite{quadrics}, we study the intersection of a general
BN-curve of degree $d$ and genus $g$ in $\pp^r$ with a hypersurface $S$ of degree $n$.
An easy dimension count (plus a tiny bit more work when $r = 2$)
implies that there are only five pairs $(r, n)$
where this intersection could be, with the exception of finitely many
$(d, g)$ pairs,
a collection of $dn$ general points on~$S$.

The main result of this paper
is that conversely, in each of these five cases, the intersection is indeed general with finitely
many exceptions:
\begin{enumerate}
\item The intersection of a plane curve with a line yields a general $d$-tuple of points on the line, always;
\item The intersection of a plane curve with a conic, always;
\item The intersection of a space curve with a quadric, with six exceptions: 
\[(d, g) \in \{(4, 1), (5, 2), (6, 2), (6, 4), (7, 5), (8, 6)\};\]
\item The intersection of a space curve with a plane, with one exception:
\[(d, g) = (6, 4)\]
\item The intersection of a curve in $\mathbb{P}^4$ with a hyperplane, with three exceptions:
\[(d, g) \in \{(8, 5), (9, 6), (10, 7)\}.\]
\end{enumerate}
In each of these exceptions, a complete description of the intersection is given.

These statements
can be reduced to statements
about the cohomology of twists of normal bundles of general BN-curves, namely that $H^1(N_C(-n)) = 0$;
as in {\bf \ref{d:aly}}, these are approached by inductive degeneration.
However, unlike in {\bf \ref{d:aly}}, we do not know of a compatible
description of the gluing data.  

The key new idea here is that when one of the curves
is contained in a hyperplane (or other hypersurface of small degree),
and certain stringent numerical constraints are satisfied,
the required properties of $N_{X \cup Y}$
can be reduced to properties of $N_X$ and $N_Y$ that (essentially)
\emph{do not depend upon the gluing data}.

\item In \cite{vogt}, Vogt shows that the normal bundle of a general BN space curve
satisfies interpolation except in exactly two cases: $(d, g) \in \{(5, 2), (6, 4)\}$.

The argument proceeds by noting that for $C$ a space curve, $H^1(N_C(-2)) = 0$ implies $N_C$ satisfies interpolation.
Using {\bf \ref{d:quadrics}}, it thus remains to show $N_C$ satisfies interpolation when $(d, g) \in \{(4, 1), (6, 2), (7, 5), (8, 6)\}$.
The cases $(d, g) \in \{(4, 1), (6, 2)\}$ are done in {\bf \ref{d:aly}}, so it
remains to show $N_C$ satisfies interpolation
when $(d, g) \in \{(7, 5), (8, 6)\}$.

In these cases, degeneration to a reducible curve is difficult, and new techniques are needed.
For curves of degree $7$ and genus $5$, which are projections of canonical curves in $\pp^4$
from a point on the curve, Vogt finds and analyzes
a description of the normal bundle
exact sequence associated to the projection,
which is compatible with the description of
a canonical curve in $\pp^4$ as the complete intersection of a net of quadrics.
For curves of degree $8$ and genus $6$, Vogt degenerates to a smooth curve lying
on a cubic surface with $3$ ordinary double points.

\item In joint work with Vogt \cite{p4}, we show that, for $C$ a general BN-curve in $\pp^4$,
the normal bundle $N_C$ (respectively the twist $N_C(-1)$)
satisfies interpolation, except in exactly $1$ case: $(d, g) = (6, 2)$ (respectively exactly $4$ cases: $(d, g) \in \{(6, 2), (8, 5), (9, 6), (10, 7)\}$).

Interpolation for $N_C(-1)$ implies that $C$ can pass through $n$ points which are general 
subject to the constraint that $d$ of them lie in a transverse hyperplane,
and subject to~\eqref{nbound}.

Unlike in $\pp^r$ for $r \geq 5$ --- where general curves
have only $1$- and $2$- secant lines --- (most) curves in $\pp^4$
have trisecant lines; the techniques of {\bf \ref{d:aly}} for inductively degenerating to
reducible curves one component
of which is a line can thus be applied here in greater generality.
Combining this with methods of {\bf \ref{d:quadrics}}, we devise an inductive argument
to prove interpolation for $N_C$ and $N_C(-1)$.

\item Finally, in \cite{ibe}, we
deduce ``bounded-error approximations'' which are valid for BN-curves of arbitrary
degree and genus, in a projective space of arbitrary dimension.

For example, we show that a BN-curve of degree $d$ and genus $g$ in $\pp^r$
passes through $n$ general points if
\[n \leq \left\lfloor \frac{(r + 1)d - (r - 3)(g - 1)}{r - 1}\right \rfloor - 3\]
(compare to \eqref{nbound}),
and gives similar statements when some of the points are constrained to lie in a hyperplane.

The proof is by degeneration to reducible BN-curves whose components
fall into a case which has already been analyzed in one of the above papers;
this degeneration is studied using methods introduced in {\bf \ref{d:quadrics}}.
\end{enumerate}

These results let us build the desired reducible curves for the maximal rank conjecture simply by showing
that there are integers (representing the degrees and genera of the components)
satisfying certain systems of inequalities --- a problem considered in Section~\ref{sec:inequalities}.

\section{The Incorrect Fiber Dimension \label{sec:fiber}}
The union described in the previous section can be constructed in the space of stable maps:
Writing $C' \cup_\Gamma C''$ for the curve obtained
from $C'$ and $C''$ by gluing along $\Gamma$,
we obtain a map $f \colon C' \cup_\Gamma C'' \to \pp^r$
(which may not be an immersion).
Conditions under which such unions are BN-curves are studied in the following
sequence of papers:

\begin{enumerate}[\bf I.]
\item \label{d:tan} In \cite{tan}, we show $T_{\pp^r}|_C$ satisfies interpolation,
where $C \subset \pp^r$ is a general BN-curve. We also give results for the twist
$T_{\pp^r}|_C(-1)$.

This implies an analog of the question considered
in Section~\ref{sec:uniform}, for maps from fixed curves with fixed marked points
which must be sent to the specified points in $\pp^r$.

As with the papers on interpolation
for normal bundles, the argument is via degeneration --- but for $T_{\pp^r}|_C$,
the gluing data is easy to understand.

\item \label{d:rbn} In \cite{rbn}, we study this construction of reducible curves $f \colon C' \cup_\Gamma C'' \to \pp^r$
in the regime where both components are nonspecial (as well as some other special cases).

First we leverage the results of {\bf \ref{d:tan}}
to calculate the fiber dimension of the map from the space of
stable maps to the moduli space of curves
at certain reducible curves, thereby showing they are BN-curves.

Then we show such reducible curves are BN-curves (subject to some mild conditions),
even when the fiber dimension is wrong,
by showing that they lie in the same component as
another curve which we know is a BN-curve by calculation of the fiber dimension.
Rather than finding an irreducible curve in the space of maps,
the key insight here is to draw
a ``broken arc'' (iteratively specialize and then deform)
in the space of stable maps, connecting these two points of the moduli space:

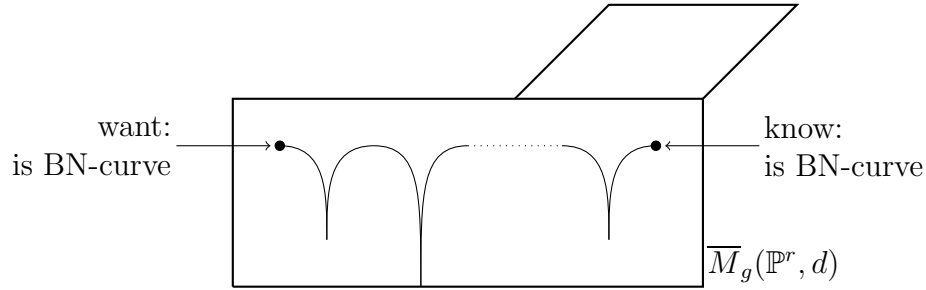
\begin{center}
\begin{tikzpicture}[scale=1.25]
\draw[thick] (0, 0) -- (0, 2) -- (5, 2) -- (5, 0) -- (0, 0);
\draw[thick] (3, 2) -- (4, 3) -- (6, 3) -- (5, 2);
\filldraw (0.5, 1.5) circle [radius=0.05];
\draw (0.5, 1.5) .. controls (1, 1.5) and (1, 1) .. (1, 0.5);
\draw (1, 0.5) .. controls (1, 1) and (1, 1.5) .. (1.5, 1.5);
\draw (1.5, 1.5) .. controls (2, 1.5) and (2, 1) .. (2, 0);
\draw (2, 0) .. controls (2, 1) and (2, 1.5) .. (2.5, 1.5);
\draw[dotted] (2.5, 1.5) -- (3.5, 1.5);
\draw (3.5, 1.5) .. controls (4, 1.5) and (4, 1) .. (4, 0.5);
\draw (4, 0.5) .. controls (4, 1) and (4, 1.5) .. (4.5, 1.5);
\filldraw (4.5, 1.5) circle [radius=0.05];
\draw[->] (-0.6, 1.5) -- (0.4, 1.5);
\draw (-1.5, 1.5) node[align=right]{want: \\ is BN-curve};
\draw[<-] (4.6, 1.5) -- (5.6, 1.5);
\draw (6.5, 1.5) node[align=left]{know: \\ is BN-curve};
\draw (5.75, 0.25) node{$\bar{M}_{g}(\pp^r, d)$};
\end{tikzpicture}
\end{center}

Provided we check the specializations are to smooth points of the space of stable maps,
this shows our given such reducible curve is in the same component as the other curve, and is thus
a BN-curve as desired.

These arcs are constructed by further specializing one of the components, say $C'$,
to a reducible curve $C'_1 \cup D_1'$; this results in a specialization of $C' \cup C''$ given by
\[(C'_1 \cup D_1') \cup C'' = C_1' \cup (D_1' \cup C'').\]
We then deform $D_1' \cup C''$ to a smooth curve $C_1''$.
Finally, we iterate this procedure, alternating between components
(next we would specialize $C_1''$ --- to a different reducible curve, not back to $D'_1 \cup C''$):

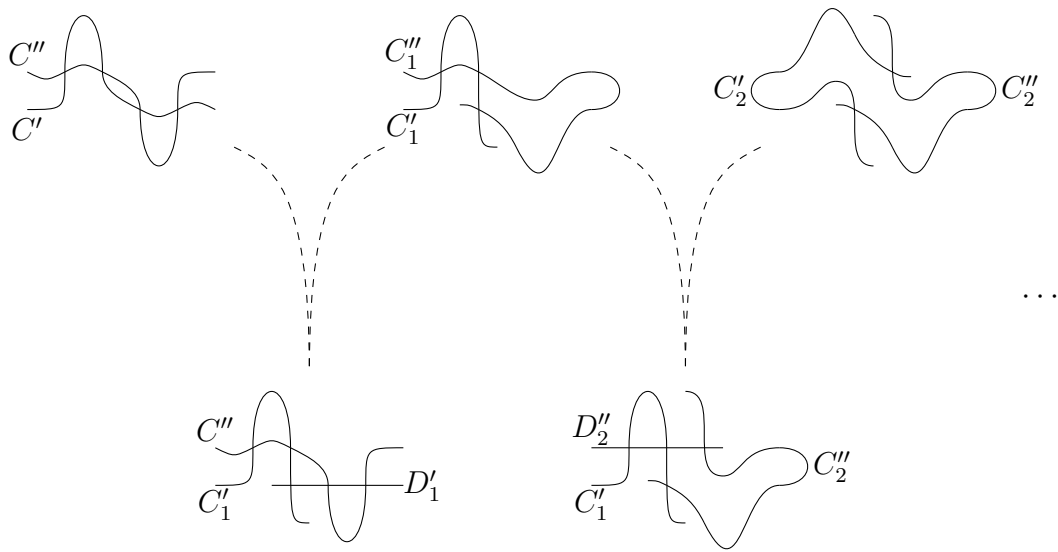
\begin{center}
\begin{tikzpicture}[scale=0.5]
\draw (0, -0.5) node{$C'$};
\draw (0, 1.5) node{$C''$};
\draw (1, 1) .. controls (1, 3) and (2, 3) .. (2, 1);
\draw (1, 1) .. controls (1, 0) .. (0, 0);
\draw (2, 1) .. controls (2, 0.5) and (2.5, 0.25) .. (3, 0);
\draw (3, 0) .. controls (3.5, -0.25) .. (4, 0);
\draw (4, 0) .. controls (4.5, 0.25) .. (5, 0);
\draw (4, 0) .. controls (4, -2) and (3, -2) .. (3, 0);
\draw (4, 0) .. controls (4, 1) .. (5, 1);
\draw (3, 0) .. controls (3, 0.5) and (2.5, 0.75) .. (2, 1);
\draw (2, 1) .. controls (1.5, 1.25) .. (1, 1);
\draw (1, 1) .. controls (0.5, 0.75) .. (0, 1);
\begin{scope}[shift={(5, -10)}]
\draw (0, -0.5) node{$C_1'$};
\draw (0, 1.5) node{$C''$};
\draw (5.5, 0) node{$D_1'$};
\draw (1, 1) .. controls (1, 3) and (2, 3) .. (2, 1);
\draw (1, 1) .. controls (1, 0) .. (0, 0);
\draw (2, 1) .. controls (2, -1) .. (2.5, -1);
\draw (1.5, 0) -- (5, 0);
\draw (4, 0) .. controls (4, -2) and (3, -2) .. (3, 0);
\draw (4, 0) .. controls (4, 1) .. (5, 1);
\draw (3, 0) .. controls (3, 0.5) and (2.5, 0.75) .. (2, 1);
\draw (2, 1) .. controls (1.5, 1.25) .. (1, 1);
\draw (1, 1) .. controls (0.5, 0.75) .. (0, 1);
\end{scope}
\begin{scope}[shift={(10, 0)}]
\draw (0, -0.5) node{$C_1'$};
\draw (0, 1.5) node{$C_1''$};
\draw (1, 1) .. controls (1, 3) and (2, 3) .. (2, 1);
\draw (1, 1) .. controls (1, 0) .. (0, 0);
\draw (2, 1) .. controls (2, -1) .. (2.5, -1);
\draw (3.5, 0.25) .. controls (3, 0.25) and (2.5, 0.75) .. (2, 1);
\draw (3.5, 0.25) .. controls (4, 0.25) and (4, 1) .. (5, 1);
\draw (5, 1) .. controls (6, 1) and (6, 0) .. (5, 0);
\draw (5, 0) .. controls (4, 0) and (4, -3) .. (3, -1);
\draw (3, -1) .. controls (2.75, -0.5) and (2.5, -0.25) .. (2, 0);
\draw (2, 0) .. controls (1.75, 0.125) .. (1.5, 0.125);
\draw (2, 1) .. controls (1.5, 1.25) .. (1, 1);
\draw (1, 1) .. controls (0.5, 0.75) .. (0, 1);
\end{scope}
\begin{scope}[shift={(15, -10)}]
\draw (0, -0.5) node{$C_1'$};
\draw (0, 1.5) node{$D_2''$};
\draw (6.4, 0.5) node{$C_2''$};
\draw (1, 1) .. controls (1, 3) and (2, 3) .. (2, 1);
\draw (1, 1) .. controls (1, 0) .. (0, 0);
\draw (2, 1) .. controls (2, -1) .. (2.5, -1);
\draw (3.5, 0.25) .. controls (2.5, 0.25) and (3.5, 2.5) .. (2.5, 2.5);
\draw (3.5, 1) -- (0, 1);
\draw (3.5, 0.25) .. controls (4, 0.25) and (4, 1) .. (5, 1);
\draw (5, 1) .. controls (6, 1) and (6, 0) .. (5, 0);
\draw (5, 0) .. controls (4, 0) and (4, -3) .. (3, -1);
\draw (3, -1) .. controls (2.75, -0.5) and (2.5, -0.25) .. (2, 0);
\draw (2, 0) .. controls (1.75, 0.125) .. (1.5, 0.125);
\end{scope}
\begin{scope}[shift={(20, 0)}]
\draw (6.4, 0.5) node{$C_2''$};
\draw (-1.3, 0.5) node{$C_2'$};
\draw (1.5, 0.75) .. controls (2.5, 0.75) and (1.5, -1.5) .. (2.5, -1.5);
\draw (1.5, 0.75) .. controls (1, 0.75) and (1, 0) .. (0, 0);
\draw (0, 0) .. controls (-1, 0) and (-1, 1) .. (0, 1);
\draw (0, 1) .. controls (1, 1) and (1, 4) .. (2, 2);
\draw (2, 2) .. controls (2.25, 1.5) and (2.5, 1.25) .. (3, 1);
\draw (3, 1) .. controls (3.25, 0.875) .. (3.5, 0.875);
\draw (3.5, 0.25) .. controls (2.5, 0.25) and (3.5, 2.5) .. (2.5, 2.5);
\draw (3.5, 0.25) .. controls (4, 0.25) and (4, 1) .. (5, 1);
\draw (5, 1) .. controls (6, 1) and (6, 0) .. (5, 0);
\draw (5, 0) .. controls (4, 0) and (4, -3) .. (3, -1);
\draw (3, -1) .. controls (2.75, -0.5) and (2.5, -0.25) .. (2, 0);
\draw (2, 0) .. controls (1.75, 0.125) .. (1.5, 0.125);
\end{scope}
\draw[dashed] (5.5, -1) .. controls (6.5, -1.5) and (7.5, -2) .. (7.5, -7);
\draw[dashed] (9.5, -1) .. controls (8.5, -1.5) and (7.5, -2) .. (7.5, -7);
\begin{scope}[shift={(10, 0)}]
\draw[dashed] (5.5, -1) .. controls (6.5, -1.5) and (7.5, -2) .. (7.5, -7);
\draw[dashed] (9.5, -1) .. controls (8.5, -1.5) and (7.5, -2) .. (7.5, -7);
\end{scope}
\draw (27, -5) node{$\cdots$};
\end{tikzpicture}
\end{center}

Note that even if 
$C'$ and $C''$ do not meet at any additional
point not in $\Gamma$, and have distinct tangent directions at the points of $\Gamma$ --- so that $f$ is the natural immersion of the scheme-theoretic union --- this broken arc may still not make sense in the Hilbert scheme
compactification, so it is important to
work in the space of stable maps even in this case.

\item Finally, in \cite{rbn2}, we repeat the analysis in {\bf \ref{d:rbn}} to study this construction of reducible curves
in the regime where such reducible curves can be constructed using
the approximate results on interpolation discussed in the Section~\ref{sec:uniform} (as well as some other special cases).
This is a separate paper from {\bf \ref{d:rbn}} since results of {\bf \ref{d:rbn}} are needed in the proof
of many of the results on interpolation discussed in Section~\ref{sec:uniform}, while those results on interpolation
are needed for \cite{rbn2}.
\end{enumerate}

\section{The Hyperplane Section \label{sec:hyp}}

Our study of subschemes of the form $C'' \cup (C' \cap S)$ which arise
in the inductive argument is divided as follows:

\begin{enumerate}
\item In \cite{hyp}, we show by degeneration that the union of hyperplane sections
\[(C_1 \cup C_2 \cup \cdots \cup C_n) \cap H\]
of independently
general BN-curves $C_1, C_2, \ldots, C_n$
satisfies maximal rank for polynomials of degree $k$, unless $k = 2$ and
some $C_i$ is special.

In low dimensions, we furthermore show that if $X \subset H$ and its hyperplane section $X \cap H'$ satisfy maximal rank,
then subject to mild conditions, so does the union $X \cup (C \cap H)$ of $X$ with the hyperplane
section of an independently general BN-curve $C$.
The proofs of these statements depend crucially on results of \cite{quadrics} discussed earlier.

\item In \cite{mrc}, one of the key steps is to study conditions under which $C'$ can be further specialized
so that its hyperplane section becomes independent from $C''$.

As an analogy, consider a set of $1$ black point and $5$ white points in the
plane, which are general subject to the condition that they lie on a conic.
The black and white points are not independent --- i.e.\ there is no description of what the white points can be that doesn't reference the position of the black point.
However, we can specialize the conic to the union of two lines, such that the black point
and $1$ white point lie on one line, while $4$ white points lie on the other line:

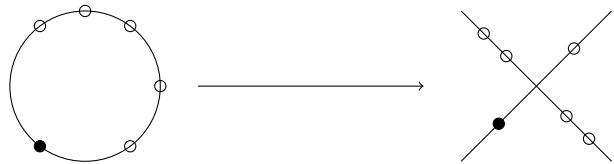
\begin{center}
\begin{tikzpicture}
\draw (0, 0) circle [radius=1];
\draw (5, -1) -- (7, 1);
\draw (5, 1) -- (7, -1);
\filldraw (-0.6, -0.8) circle[radius=0.075];
\draw (0.6, 0.8) circle[radius=0.075];
\draw (-0.6, 0.8) circle[radius=0.075];
\draw (0.6, -0.8) circle[radius=0.075];
\draw (1, 0) circle[radius=0.075];
\draw (0, 1) circle[radius=0.075];
\filldraw (5.5, -0.5) circle[radius=0.075];
\draw (6.5, 0.5) circle[radius=0.075];
\draw (5.3, 0.7) circle[radius=0.075];
\draw (5.6, 0.4) circle[radius=0.075];
\draw (6.4, -0.4) circle[radius=0.075];
\draw (6.7, -0.7) circle[radius=0.075];
\draw[->] (1.5, 0) -- (4.5, 0);
\end{tikzpicture}
\end{center}

After specialization, the black and white points become independent:
The white points specialize to a set of $5$ points which are general subject to the constraint
that $4$ of them are collinear --- a description that doesn't reference the position of the black point.

In our setting, we further specialize $C'$ to a reducible curve $C_1' \cup C_2'$,
such that
\begin{equation} \label{ind}
[C_2' \cap H] \times [C_1' \cap C_2'] \in \sym^{\deg C_2'} \pp^r \times \sym^{\# (C_1' \cap C_2')} H
\end{equation}
is general. This induces a specialization of $C' \cup C''$:

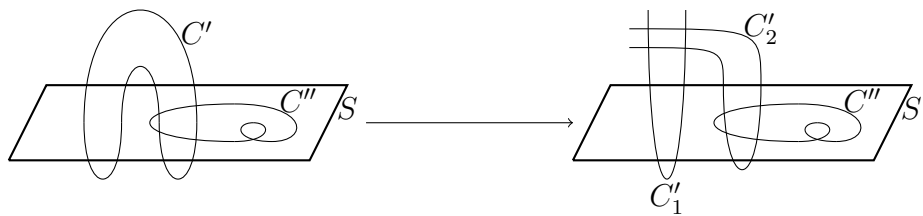
\begin{center}
\begin{tikzpicture}[scale=0.5]
\draw[thick] (0, 1) -- (1, 3) -- (9, 3) -- (8, 1) -- (0, 1);
\draw (2, 2) .. controls (2, 6) and (5, 6) .. (5, 2);
\draw (2, 2) .. controls (2, 0) and (3, 0) .. (3, 2);
\draw (4, 2) .. controls (4, 0) and (5, 0) .. (5, 2);
\draw (3, 2) .. controls (3, 4) and (4, 4) .. (4, 2);
\draw (6, 2.5) .. controls (3, 2.5) and (3, 1.5) .. (6, 1.5);
\draw (6, 2.5) .. controls (8, 2.5) and (8, 1.5) .. (7, 1.5);
\draw (6, 1.5) .. controls (7, 1.5) and (7, 2) .. (6.5, 2);
\draw (7, 1.5) .. controls (6, 1.5) and (6, 2) .. (6.5, 2);
\draw (9, 2.4) node{$S$};
\draw (5, 4.4) node{$C'$};
\draw (7.7, 2.6) node{$C''$};
\begin{scope}[shift={(15, 0)}]
\draw[thick] (0, 1) -- (1, 3) -- (9, 3) -- (8, 1) -- (0, 1);
\draw (2, 5) .. controls (2, -1) and (3, -1) .. (3, 5);
\draw (4, 3) .. controls (4, 0) and (5, 0) .. (5, 3);
\draw (4, 3) .. controls (4, 4) .. (1.5, 4);
\draw (5, 3) .. controls (5, 4.5) .. (1.5, 4.5);
\draw (6, 2.5) .. controls (3, 2.5) and (3, 1.5) .. (6, 1.5);
\draw (6, 2.5) .. controls (8, 2.5) and (8, 1.5) .. (7, 1.5);
\draw (6, 1.5) .. controls (7, 1.5) and (7, 2) .. (6.5, 2);
\draw (7, 1.5) .. controls (6, 1.5) and (6, 2) .. (6.5, 2);
\draw (9, 2.4) node{$S$};
\draw (5, 4.5) node{$C_2'$};
\draw (7.7, 2.6) node{$C''$};
\draw (2.5, 0) node{$C_1'$};
\end{scope}
\draw[->] (9.5, 2) -- (15, 2);
\end{tikzpicture}
\end{center}

Via the method of Hirschowitz, we reduce maximal rank for $C' \cup C''$
to maximal rank for $C'$, and maximal rank for $C'' \cup (C' \cap H)$,
which in turn reduces to maximal rank for $C'' \cup (C_1' \cap H) \cup (C_2' \cap H)$.
Under our assumption \eqref{ind}, $C''$, $C_1' \cap H$, and $C_2' \cap H$
are an \emph{independently} general BN-curve, hyperplane section of a BN-curve, and set of points.

The upshot is that we can then argue by induction on
the following stronger hypothesis.
(Note that taking $n = \epsilon = 0$
recovers the maximal rank conjecture.)

\begin{thm} \label{main}
Fix an inclusion $\pp^r \subset \pp^{r + 1}$ (for $r \geq 3$), and let $k$ be a positive integer.
Let $C \subset \pp^r$ be a general BN-curve or a general degenerate rational curve
of degree $0 < d < r$.
Let $D_1, D_2, \ldots, D_n \subset \pp^{r + 1}$ be independently
general BN-curves, which are required to be nonspecial if $k = 2$ and $r \geq 4$.
Let $p_1, p_2, \ldots, p_\epsilon \in \pp^r$
be a general set of points.
Then any subset of
\[T \colonequals C \cup ((D_1 \cup D_2 \cup \cdots \cup D_n) \cap \pp^r) \cup \{p_1, p_2, \ldots, p_\epsilon\} \subset \pp^r\]
which contains $C$
satisfies maximal rank for polynomials of degree $k$.
\end{thm}
\end{enumerate}

\section{\label{sec:inequalities} Integer Solutions to Systems of Inequalities}

Applying the techniques of previous sections, we show in \cite{mrc} that the maximal rank conjecture
may be reduced to several instances of the following problem:
\emph{Given integers (e.g.\ $r, k, d, g, \ldots$) satisfying a certain system of inequalities
(e.g.\ $\rho(d, g, r) \geq 0, \ldots$), show that there are either additional integers
(e.g.\ $d', g', d'', g'', \ldots$)
satisfying a first additional system of inequalities, or other additional integers
(e.g.\ $d_1', g_1', d_2', g_2', d'', g'', \ldots$)
satisfying a second additional system of inequalities.}

Crucially, all the inequalities arising in the proof of the maximal rank conjecture
are linear in all the variables except $r$ and $k$, with coefficients that are polynomials
in $r$, $k$, and the binomial coefficient $\binom{r + k}{k}$. For fixed $r$ and $k$,
these systems of inequalities describe compact convex polyhedra.

These computations can be approached in three steps:

\begin{enumerate}
\item First, we eliminate the additional variables one by one, using the following fact:
There exists a real number, respectively integer, $n$ satisfying the inequalities
\[n \leq \frac{a_i}{b_i} \tand n \geq \frac{c_j}{d_j}\]
if, for each $(i, j)$,
\[a_i d_j - b_i c_j \geq 0 \quad \text{respectively} \quad a_i d_j - b_i c_j \geq (b_i - 1)(d_j - 1).\]

\item Then we reduce the given problem to checking positivity of polynomials
in $r$, $k$, and $\binom{r + k}{k}$, using the following fact:
Let $P \subset \rr^n$ be a compact convex polyhedron, and $C_1, C_2 \subset \rr^n$ be
convex sets. Then $P \subset C_1 \cup C_2$ if and only if:
\begin{enumerate}
\item \label{vertex} Every vertex of $P$ is contained in either $C_1$ or $C_2$; and
\item \label{edge} Every edge of $P$ joining a vertex not contained in $C_1$ to a vertex not contained in $C_2$
meets $C_1 \cap C_2$.
\end{enumerate}

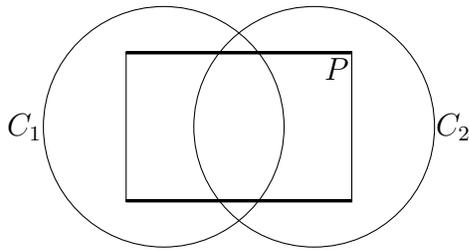
\begin{center}
\begin{tikzpicture}
\draw (0.5, 1) circle[radius=1.6];
\draw (2.5, 1) circle[radius=1.6];
\draw (0, 0) -- (0, 2);
\filldraw (0, 2) -- (3, 2) -- (3, 1.97) -- (0, 1.97) -- (0, 2);
\filldraw (0, 0) -- (3, 0) -- (3, 0.03) -- (0, 0.03) -- (0, 0);
\draw (3, 2) -- (3, 0);
\draw (-1.35, 1) node{$C_1$};
\draw (4.35, 1) node{$C_2$};
\draw (2.8, 1.77) node{$P$};
\end{tikzpicture}
\end{center}

\item Finally, we verify positivity of these polynomials, using the following fact:
A polynomial $P(r, k)$
in two variables $r$ and $k$ is positive for all $r \geq r_0$
and $k \geq k_0$, provided that:
\begin{enumerate}
\item \label{epos} Every monomial on the outside of the Newton polygon has positive
coefficient.

\item \label{sposk} The leading coefficient with respect to $r$ is positive for $k \geq k_0$.

\item \label{sposr} The leading coefficient with respect to $k$ is positive for $r \geq r_0$.

\item \label{boundary} The polynomial $P(r_0, k)$ is positive for
$k \geq k_0$.

\item \label{branch} The value of $k_0$ exceeds all branch points of the projection
onto the $k$-axis of $P(r, k) = 0$.
\end{enumerate}
\end{enumerate}

Sometimes this method may fail. For example, when $r = 17$ and $k = 4$,
there is a vertex of a polyhedron $P$ appearing in the second
step --- corresponding to the value of $d$ and $g$ for which
$\rho(d, g, 17) = 0$ and the maximal rank map is expected to be an isomorphism ---
which is not contained in either convex set. However, this vertex has non-integer values of
$d$ and $g$, and is only barely not contained in either convex set; brute force
search shows that, in this case, the desired existence of additional integers
holds for every integral $(d, g)$.
(In particular, we see that this proof of the maximal rank conjecture barely works; with
only slightly worse approximate results on interpolation, it would fail.)

Proofs of the above facts, as well as
computer code implementing this method (combined with brute-force search
where it fails), are given in Appendix~E of~\cite{mrc}.

\subsection*{Acknowledgements}

The author would like to thank Joe Harris for
his guidance throughout this research,
as well as Atanas Atanasov, Edoardo Ballico, Brian Osserman, Sam Payne, Ravi Vakil, Isabel Vogt, David Yang, and other members of the Harvard and MIT mathematics departments,
for helpful conversations or comments on this manuscript.
The author would also like
to acknowledge the generous
support both of the Fannie and John Hertz Foundation,
and of the Department of Defense
(DoD) through the National  Defense Science and Engineering Graduate Fellowship (NDSEG) Program.

\bibliographystyle{amsplain.bst}
\bibliography{mrcbib}

\providecommand{\bysame}{\leavevmode\hbox to3em{\hrulefill}\thinspace}
\providecommand{\MR}{\relax\ifhmode\unskip\space\fi MR }
\providecommand{\MRhref}[2]{%
  \href{http://www.ams.org/mathscinet-getitem?mr=#1}{#2}
}
\providecommand{\href}[2]{#2}
\begin{thebibliography}{10}

\bibitem{aly}
Atanas Atanasov, Eric Larson, and David Yang, \emph{Interpolation for normal
  bundles of general curves}, accepted for publication in \emph{Memoirs of the
  AMS}; preprint available at \url{http://arxiv.org/abs/1509.01724v2}.

\bibitem{gp}
David Gieseker, \emph{Stable curves and special divisors: {P}etri's
  conjecture}, Invent. Math. \textbf{66} (1982), no.~2, 251--275. \MR{656623
  (83i:14024)}

\bibitem{bn}
Phillip Griffiths and Joseph Harris, \emph{On the variety of special linear
  systems on a general algebraic curve}, Duke Math. J. \textbf{47} (1980),
  no.~1, 233--272. \MR{563378 (81e:14033)}

\bibitem{severi}
Joeseph Harris, \emph{On the {S}everi problem}, Invent. Math. \textbf{84}
  (1986), no.~3, 445--461. \MR{837522}

\bibitem{hh}
R.~Hartshorne and A.~Hirschowitz, \emph{Smoothing algebraic space curves},
  Algebraic geometry, {S}itges ({B}arcelona), 1983, Lecture Notes in Math.,
  vol. 1124, Springer, Berlin, 1985, pp.~98--131. \MR{805332 (87h:14023)}

\bibitem{kl}
Steven Kleiman and Dan Laksov, \emph{On the existence of special divisors},
  Amer. J. Math. \textbf{94} (1972), 431--436. \MR{0323792}

\bibitem{quadrics}
Eric Larson, \emph{The generality of a section of a curve},
  \url{http://arxiv.org/abs/1605.06185}.

\bibitem{ibe}
\bysame, \emph{Interpolation with bounded error},
  \url{https://arxiv.org/abs/1711.01729}.

\bibitem{tan}
\bysame, \emph{Interpolation for restricted tangent bundles of general curves},
  Algebra Number Theory \textbf{10} (2016), no.~4, 931--938. \MR{3519101}

\bibitem{rbn}
\bysame, \emph{Constructing reducible {B}rill--{N}oether curves},
  \url{https://arxiv.org/abs/1603.02301}.

\bibitem{rbn2}
\bysame, \emph{Constructing reducible {B}rill--{N}oether curves {II}},
  \url{https://arxiv.org/abs/1711.02752}.

\bibitem{hyp}
\bysame, \emph{The maximal rank conjecture for sections of curves},
  \url{https://arxiv.org/abs/1208.2730}.

\bibitem{mrc}
\bysame, \emph{The maximal rank conjecture},
  \url{https://arxiv.org/abs/1711.04906}.

\bibitem{p4}
Eric Larson and Isabel Vogt, \emph{Interpolation for curves in $\pp^4$},
  \url{https://arxiv.org/abs/1708.00028}.

\bibitem{sernesi}
Edoardo Sernesi, \emph{On the existence of certain families of curves}, Invent.
  Math. \textbf{75} (1984), no.~1, 25--57. \MR{728137}

\bibitem{vogt}
Isabel Vogt, \emph{Interpolation for {B}rill--{N}oether space curves}, To
  appear in \emph{Manuscripta Mathematica}; preprint available at
  \url{https://arxiv.org/abs/1611.00081}.

\end{thebibliography}

\end{document}